\documentclass [a4 paper,12pt]{article}
\usepackage{mathrsfs}

\usepackage[usenames]{color}
\usepackage{bbm,amsmath,amssymb,amsfonts,amsthm,graphics,graphicx}

\newtheorem{lem}{Lemma}
\newtheorem{thm}[lem]{Theorem}

\newtheorem{con}{Conjecture}

\DeclareMathOperator{\Pro}{\mathbbm{P}}
\DeclareMathOperator{\V}{\mathbbm{Var}}
\DeclareMathOperator{\E}{\mathbbm{E}}
\DeclareMathOperator{\En}{\mathscr{E}}
\DeclareMathOperator{\EnL}{\mathscr{E}_{\mathit L}}
\DeclareMathOperator{\La}{\overline{\mathbf{L}}}
\DeclareMathOperator{\Lf}{\mathbf{L}_1}
\DeclareMathOperator{\Ls}{\mathbf{L}_2}

\DeclareMathOperator{\si}{\sigma}
\DeclareMathOperator{\la}{\lambda}

\title{The Laplacian energy of random graphs\footnote{Supported by NSFC No.10831001, PCSIRT and the
``973" program. }}
\author{
\small Wenxue Du, Xueliang Li, Yiyang Li\\
\small Center for Combinatorics and LPMC-TJKLC \\
\small Nankai University, Tianjin 300071, China \\
\small Email: lxl@nankai.edu.cn}
\date{}

\begin{document}
\maketitle

\begin{abstract}
Gutman {\it et al.} introduced the concepts of energy $\En(G)$
and Laplacian energy $\EnL(G)$ for a simple graph $G$, and
furthermore, they proposed a conjecture that for every graph $G$,
$\En(G)$ is not more than $\EnL(G)$. Unfortunately, the
conjecture turns out to be incorrect since Liu {\it et al.} and
Stevanovi\'{c} {\it et al.} constructed counterexamples. However,
So {\it et al.} verified the conjecture for bipartite graphs. In
the present paper, we obtain, for a random graph, the lower and
upper bounds of the Laplacian energy, and show that the
conjecture is true for almost all graphs.\\[3mm]

\end{abstract}

\section{Introduction}

Throughout this paper, $G$ denotes a simple graph of order $n$. The
eigenvalues $\lambda_1,\ldots,\lambda_n$ of the adjacency matrix
$\mathbf{A}(G)=(a_{ij})_{n\times n}$ are said to be the {\em
eigenvalues of $G$}. In chemistry, there is a closed relation
between the molecular orbital energy levels of $\pi$-electrons in
conjugated hydrocarbons and the eigenvalues of the corresponding
molecular graph. For the H\"{u}chkel molecular orbital
approximation, the total $\pi$-electron energy in conjugated
hydrocarbons is given by the sum of absolute values of the
eigenvalues corresponding to the molecular graph $G$ in which the
maximum degree is not more than 4 in general. In 1970s, Gutman
\cite{G} extended the concept of energy to all simple graphs $G$,
and defined that
$$\En(G)=\sum_{i=0}^n|\lambda_i|,$$
where $\lambda_1,\ldots,\lambda_n$ are the eigenvalues of $G$.
Evidently, one can immediately get the energy of a graph by
computing the eigenvalues of the graph. It is rather hard, however,
to compute the eigenvalues for a large matrix, even for a large
symmetric (0,1)-matrix like $\mathbf{A}(G)$. So many researchers
established a lot of lower and upper bounds to estimate the
invariant for some classes of graphs. For further details, we refer
readers to the comprehensive survey \cite{glz}. But there is a
common flaw for those inequalities that only a few graphs attain the
equalities of those bounds. Consequently we can hardly see the major
behavior of the invariant $\En(G)$ for most graphs with respect to
other graph parameters ($|V(G)|$, for instance). In the next
section, however, we shall present an exact estimate of the energy
for almost all graphs by Wigner's semi-circle law.

In spectral graph theory, the matrix
$\mathbf{L}(G)=\mathbf{D}(G)-\mathbf{A}(G)$ is called  {\em
Laplacian matrix of $G$}, where $\mathbf{D}(G)$ is a diagonal matrix
in which $d_{ii}$ equals the degree $d_G(v_i)$ of the vertex $v_i$,
$i=1,\ldots,n$. Gutman {\it et al.} \cite{GZ} introduced a new
matrix $\La(G)$ for a simple graph $G$, {\it i.e.,}
$$\La(G)=\mathbf{L}(G)-\sum_{i=1}^nd_G(v_i)/n\mathbf{I}_n=
\mathbf{L}(G)-2\sum_{i=1}^n\sum_{j>i}a_{ij}/n\mathbf{I}_n,$$
where $\mathbf{I}_n$ is the unit matrix of order $n$, and defined
the {\em Laplacian energy $\EnL(G)$ of $G$}, {\it i.e.,}
$$\EnL(G)=\sum_{i=1}^n|\zeta_i|,$$
where $\zeta_1,\ldots,\zeta_n$ are the eigenvalues of $\La(G)$.
Obviously, we can easily evaluate the Laplacian energy $\EnL(G)$ if
we could obtain the eigenvalues of $\La(G)$. In Section 3 we shall
establish the lower and upper bounds of the Laplacian energy for
almost all graphs by exploring the spectral distribution of the
matrix $\La(G_n(p))$ for a random graph $G_n(p)$ constructed from
the classical Erd\"{o}s--R\'{e}nyi model (see \cite{BB}).

In a recent paper \cite{gut}, Gutman {\it et al.} proposed the
following conjecture concerning the relation between the energy and
the Laplacian energy of a graph.

\begin{con}\label{Gutman} Let $G$ be a simple graph. Then
$\En(G)\le\EnL(G).$
\end{con}

Unfortunately, the conjecture turns out to be incorrect. In fact,
Liu {\it et al.} \cite{LL} and Stevanovi\'{c} {\it et al.}
\cite{SSM} constructed two classes of graphs violating the
assertion. However, So {\it et al.} \cite{sw} proved that the
conjecture is true for bipartite graphs. We shall show that the
conjecture above is true for almost all graphs by comparing the
energy with the Laplacian energy of a random graph in the third
section.

\section{The energy of $G_n(p)$}

 In this section, we shall formulate an exact estimate of the
 energy for almost all graphs by Wigner's semi-circle law.

 We start by recalling the Erd\"{o}s--R\'{e}nyi
 model $\mathcal{G}_{n}(p)$ (see \cite{BB}), which consists of
 all graphs with vertex set $[n]=\{1,2,\ldots,n\}$ in which the
 edges are chosen independently with probability $p=p(n)$.
 Apparently, the adjacency matrix $\mathbf{A}(G_n(p))$
 of the random graph $G_n(p)\in\mathcal{G}_{n}(p)$ is a random
 matrix, and thus one can readily evaluate the energy of
 $G_n(p)$ once the spectral distribution of the random matrix
  $\mathbf{A}(G_n(p))$ is
 known.

 In fact, the research on the spectral distributions of random
 matrices is rather abundant and active, which can be traced
 back to \cite{Wi}. We refer readers to \cite{B, De, Me} for an overview and some
 spectacular progress in this field. One important
 achievement in that field is Wigner's semi-circle law which
 characterizes the limiting spectral distribution of the empirical
 spectral distribution of eigenvalues for a sort of random
 matrix.

 In order to characterize the statistical properties of the wave
 functions of quantum mechanical systems, Wigner
 in 1950s investigated the spectral distribution for a sort of random
 matrix, so-called {\em Wigner matrix},
 $$\mathbf{X}_n:=(x_{ij}), ~~1\le i,j\le n,$$
 which satisfies the
following properties:\begin{itemize}

  \item $x_{ij}$'s are independent random variables with $x_{ij}=x_{ji}$;

  \item the $x_{ii}$'s have the same distribution $F_1$,
  while the $x_{ij}$'s $(i\neq j)$ are to
 possess
 the same distribution $F_2$;

 \item $\V(x_{ij})=\si_2^2<\infty$ for all $1\le i< j\le
 n$.
 \end{itemize}
 We denote the eigenvalues of
 $\mathbf{X}_n$ by $\la_{1,n},\la_{2,n},\ldots,\la_{n,n}$, and their
 empirical spectral distribution (ESD) by
 $$\Phi_{\mathbf{X}_n}(x)=
 \frac{1}{n}\cdot\#\{\la_{i,n}\mid\la_{i,n}\le x, ~i=1,2,\ldots, n\}.$$
 Wigner \cite{W55,W58} considered the limiting spectral distribution (LSD)
 of $\mathbf{X}_n$, and obtained the semi-circle law.
\begin{thm}\label{Thm-1} Let $\mathbf{X}_n$
be a Wigner matrix. Then
 $$\lim_{n\rightarrow\infty}\Phi_{n^{-1/2}\mathbf{X}_n}(x)=\Phi(x)
\mbox{ a.s. }$${\it i.e.,} with probability 1, the ESD
$\Phi_{n^{-1/2}\mathbf{X}_n}(x)$ converges weakly to a distribution
$\Phi(x)$ as $n$ tends to infinity, where $\Phi(x)$ has the density
 $$\phi(x)=
  \frac{1}{2\pi\si_2^2}\sqrt{4\si_2^2- x^2 }~\mathbf{1}_{|x|\le2\si_2}.$$
\end{thm}

\noindent{\bf Remark.} It is interesting that the existence of the
second moment of the off-diagonal entries is the necessary and
sufficient condition for the semi-circle law, and there is no moment
requirement on the diagonal elements. Furthermore, we can get more
information about spectra of Wigner matrices. Set $\mu_i=\int
x\hspace{2pt}dF_i$ $(i=1,2)$ and
$$\overline{\mathbf{X}}_n
=\mathbf{X}_n-\mu_1\mathbf{I}_n-\mu_2(\mathbf{J}_n-\mathbf{I}_n),$$
where $\mathbf{J}_n$ is the all 1's matrix. One can easily check
that each entry of $\overline{\mathbf{X}}_n$ has mean 0. By means of
Wigner's trace method, one can show that the spectral radius
$\rho(n^{-1/2}\overline{\mathbf{X}}_n)$ converges to $2\sigma_2$
with probability 1 as $n$ tends to infinity (see Theorem 2 in
\cite{FK}, for instance). For further comments on Wigner's
semi-circle law, we refer readers to the extraordinary survey by Bai
\cite{B}.\vspace{6pt}

Following the book \cite{BB}, we will say that {\it almost every}
(a.e.) graph in $\mathcal{G}_n(p)$ has a certain property $Q$ if the
probability that a random graph $G_n(p)$ has the property $Q$
converges to 1 as $n$ tends to infinity. Occasionally, we shall
write {\it almost all} instead of almost every. It is easy to see
that if $F_1$ is a {\em pointmass at 0}, {\it i.e.,} $F_1(x)=1$ for
$x\ge 0$ and $F_1(x)=0$ for $x<0$, and $F_2$ is the {\em Bernoulli
distribution with mean $p$}, then the Wigner matrix $\mathbf{X}_n$
coincides with the adjacency matrix of $G_n(p)$. Obviously,
$\si_2=\sqrt{p(1-p)}$ in this case.

To establish the exact estimate of the energy $\En(G_n(p))$ for a.e.
graph $G_n(p)$, we first present some notions and assertions. In
what follows, we shall use $\mathbf{A}$ to denote the adjacency
matrix $\mathbf{A}(G_n(p))$ for convenience. Set
$$\overline{\mathbf{A}}=\mathbf{A}-p(\mathbf{J}_n-\mathbf{I}_n).$$
Evidently, $\overline{\mathbf{A}}$ is a Wigner matrix. By means of
Theorem \ref{Thm-1}, we have \begin{equation}\label{Equ-6}
\lim_{n\rightarrow\infty}\Phi_{n^{-1/2}\overline{\mathbf{A}}}(x)=\Phi(x)\mbox{
a.s.}\end{equation} It is easy to check that each entry of
$\overline{\mathbf{A}}$ has mean 0. According to the remark above,
\begin{equation}\label{SpecRad}
\lim_{n\rightarrow\infty}\rho(n^{-1/2}\overline{\mathbf{A}})
=2\sigma_2\mbox{ a.s. }\end{equation}

We further define the {\em energy $\En(\mathbf{M})$ of a matrix
$\mathbf{M}$} as the sum of absolute values of the eigenvalues of
$\mathbf{M}$. By virtue of Equation (\ref{Equ-6}) and
(\ref{SpecRad}), we shall formulate an estimate of the energy
$\En(\overline{\mathbf{A}})$, and then establish the exact estimate
of $\En(\mathbf{A})=\En(G_n(p))$ using Lemma \ref{ky}.

According to Equation (\ref{SpecRad}), for any given $\epsilon>0$,
there exists an integer $N$ such that with probability 1, for all
$n>N$ the spectral radius $\rho(n^{-1/2}\overline{\mathbf{A}})$ is
not more than $2\sigma_2+\epsilon$. Since the density $\phi(x)$ of
$\Phi(x)$ is bounded on $\mathbb{R}$, invoking Equation
(\ref{Equ-6}) and bounded convergence theorem yields that for all
$n>N$,
\begin{eqnarray}
\lim_{n\rightarrow\infty} \int
|x|d\Phi_{n^{-1/2}\overline{\mathbf{A}}}(x)
&=&\lim_{n\rightarrow\infty}
\int_{-2\sigma_2-\epsilon}^{2\sigma_2+\epsilon}
|x|d\Phi_{n^{-1/2}\overline{\mathbf{A}}}(x)\mbox{
a.s.}\nonumber\\
&=& \int_{-2\sigma_2-\epsilon}^{2\sigma_2+\epsilon}
|x|d\Phi(x)\mbox{ a.s. }
\label{Convergence} \\
&=& \int |x|d\Phi(x).\nonumber\end{eqnarray}

We now turn to the estimate of the energy
$\En(\overline{\mathbf{A}})$. Suppose
$\overline{\lambda}_1,\ldots,\overline{\lambda}_n$ and
$\overline{\lambda}_1',\ldots,\overline{\lambda}_n'$ are the
eigenvalues of $\overline{\mathbf{A}}$ and
$n^{-1/2}\overline{\mathbf{A}}$, respectively. Clearly,
$\sum_{i=1}^n|\overline{\lambda}_i|=n^{1/2}\sum_{i=1}^n|\overline{\lambda}_i'|$.
By Equation (\ref{Convergence}), we can deduce that
\begin{eqnarray*}
  \En\left(\overline{\mathbf{A}}\right)/n^{3/2} &=&
  \frac{1}{n^{3/2}}\sum_{i=1}^n|\overline{\lambda}_i|\\
  &=&
  \frac{1}{n}\sum_{i=1}^n|\overline{\lambda}_i'|\\
  &=& \int |x|d\Phi_{n^{-1/2}\overline{\mathbf{A}}}(x)\\
  &\rightarrow& \int
  |x|d\Phi(x)\mbox{ a.s. } (n\rightarrow\infty)\\
  &=&\frac{1}{2\pi\sigma_2^2}\int_{-2\sigma_2}^{2\sigma_2}
  |x|\sqrt{4\sigma_2^2-x^2}~dx\\
  &=&
  \frac{8}{3\pi}\sigma_2
  =\frac{8}{3\pi}\sqrt{p(1-p)}.
 \end{eqnarray*}
Therefore, with probability 1, the energy
$\En\left(\overline{\mathbf{A}}\right)$ enjoys the equation as
follows:
$$\En\left(\overline{\mathbf{A}}\right)
=n^{3/2}\left(\frac{8}{3\pi}\sqrt{p(1-p)}+o(1)\right).$$

We proceed to investigate $\En(\mathbf{A})=\En(G_n(p))$ and present
the following result due to Fan.
\begin{lem}[\bf Fan \cite{ky}]\label{ky}
Let $\mathbf{X},\mathbf{Y},\mathbf{Z}$ be real symmetric matrices of
order $n$ such that $\mathbf{X}+\mathbf{Y}=\mathbf{Z}$, then
$$\sum_{i=1}^n |\lambda_i(\mathbf{X})|+\sum_{i=1}^n
|\lambda_i(\mathbf{Y})| \geq\sum_{i=1}^n |\lambda_i(\mathbf{Z})|$$
where $\lambda_i(\mathbf{M})$ $(i=1,\cdots,n)$ is an eigenvalue of
the matrix $\mathbf{M}$.
\end{lem}
It is not difficult to verify that the eigenvalues of the matrix
$\mathbf{J}_n-\mathbf{I}_n$ are $n-1$ and $-1$ of $n-1$ times.
Consequently $\En(\mathbf{J}_n-\mathbf{I}_n)=2(n-1)$. One can
readily see that
$\En\big(p(\mathbf{J}_n-\mathbf{I}_n)\big)=p\En(\mathbf{J}_n-\mathbf{I}_n)$.
Thus,
$$\En\big(p(\mathbf{J}_n-\mathbf{I}_n)\big)=2p(n-1).$$ Since
$\mathbf{A}=\overline{\mathbf{A}}+p(\mathbf{J}_n-\mathbf{I}_n)$, it
follows from Lemma \ref{ky} that with probability 1,
\begin{eqnarray*}\En(\mathbf{A})
&\le&\En\left(\overline{\mathbf{A}}\right)+\En(p(\mathbf{J}_n-\mathbf{I}_n))\\
&=&n^{3/2}\left(\frac{8}{3\pi}\sqrt{p(1-p)}+o(1)\right)+2p(n-1).
\end{eqnarray*}Consequently,
\begin{equation}\label{Equ-4}\lim_{n\rightarrow\infty}\En(\mathbf{A})/n^{3/2}
\le\frac{8}{3\pi}\sqrt{p(1-p)}\mbox{ a.s.}\end{equation} On the
other hand, since
$\overline{\mathbf{A}}=\mathbf{A}+p\big(-(\mathbf{J}_n-\mathbf{I}_n)\big)$,
we can deduce by Lemma \ref{ky} that with probability 1,
\begin{eqnarray*}\En(\mathbf{A})
&\ge&\En\left(\overline{\mathbf{A}}\right)
-\En(p\big(-(\mathbf{J}_n-\mathbf{I}_n)\big))\\
&=&\En\left(\overline{\mathbf{A}}\right)-\En(p(\mathbf{J}_n-\mathbf{I}_n))\\
&=&n^{3/2}\left(\frac{8}{3\pi}\sqrt{p(1-p)}+o(1)\right)-2p(n-1).
\end{eqnarray*}
Consequently,
\begin{equation}\label{Equ-5}\lim_{n\rightarrow\infty}\En(\mathbf{A})/n^{3/2}
\ge\frac{8}{3\pi}\sqrt{p(1-p)}\mbox{ a.s.}\end{equation} Combining
equations (\ref{Equ-4}) with (\ref{Equ-5}), we have
$$\En(\mathbf{A})
=n^{3/2}\left(\frac{8}{3\pi}\sqrt{p(1-p)}+o(1)\right)\mbox{ a.s.}$$
Recalling that $\mathbf{A}$ is the adjacency matrix of $G_n(p)$, we
thus obtain that a.e. random graph $G_n(p)$ enjoys the equation as
follows:
$$\En(G_n(p))
=n^{3/2}\left(\frac{8}{3\pi}\sqrt{p(1-p)}+o(1)\right).$$ Note that
for $p=\frac 1 2$, Nikiforov in \cite{N} got the above equation.
Here, our result is for any probability $p$, which could be seen as
a generalization of his result.

\section{The Laplacian energy of $G_n(p)$}

In this section, we shall establish the lower and upper bounds of
the Laplacian energy of $G_n(p)$ by employing the LSD of Markov
matrix. Finally, we shall show that Conjecture \ref{Gutman} is true
for almost all graphs by comparing the energy with the Laplacian
energy of a random graph.

\subsection{The limiting spectral distribution}

We begin with another random matrix we are interested in. Define a
random matrix $\mathbf{M}_n=\mathbf{X}_n-\mathbf{D}_n$ to be a {\em
Markov matrix} if $\mathbf{X}_n$ is a Wigner matrix such that $F_1$
is the pointmass at zero, and $\mathbf{D}_n$ is a diagonal matrix in
which $d_{ii}=\sum_{j\neq i}x_{ij},~i=1,\ldots,n.$ The matrix is
introduced as the derivative of a transition matrix in a Markov
process.  Bryc {\it et al.} in \cite{BDJ} obtained the LSD of Markov
matrix.  Define the standard semi-circle distribution
$\Phi_{0,1}(x)$ of zero mean and unit variance to be the measure on
the real set of compact support with density
$\phi_{0,1}(x)=\frac{1}{2\pi}\sqrt{4-x^2}~\mathbf{1}_{|x|\le 2}$.
\begin{thm}[\bf Bryc {\it et al.} \cite{BDJ}]\label{LSD-Markov}
Let $\mathbf{M}_n$ be a markov matrix such that $\int x dF_2(x)=0$
and $\sigma_2=1$. Then
$$\lim_{n\rightarrow\infty}\Phi_{n^{-1/2}\mathbf{M}_n}(x)
=\Psi(x)\mbox{ a.s.}$$ where $\Psi(x)$ is the free convolution of
the standard semi-circle distribution $\Phi_{0,1}(x)$ and the
standard normal measures. Moreover, this measure $\Psi(x)$ is a
non-random symmetric probability measure with smooth bounded
density, and does not depend on the distribution of the random
variable $x_{ij}$.
\end{thm}
\noindent{\bf Remark.} To prove the theorem above, Bryc {\it et al.}
employ the moment approach. In fact, they show that for each
positive integer $k$,
\begin{equation}\label{Moments}
\lim_{n\rightarrow\infty}\int
x^k~d\Phi_{n^{-1/2}\mathbf{M}_n}(x)=\int x^k~d\Psi(x)\mbox{ a.s.
}\end{equation}  For two probability measures $\mu$ and $\nu$, there
exists a unique probability measure $\mu\boxplus\nu$ called the {\em
free convolution} of $\mu$ and $\nu$. This concept introduced by
Voiculescu \cite{voi} via $C^*$-algebraic will be discussed in
detail in the second part of this section.

Let $G_n(p)$ be a random graph of $\mathcal{G}_{n}(p)$. Set
$\sigma=\sqrt{p(1-p)}$. One can easily see that $\sigma^2$ is the
variance of the random variable $a_{ij}$ $(i>j)$ in
$\mathbf{A}(G_n(p))$. To state the main result of this part, we
present a new matrix $\Lf$ as follows:
\begin{eqnarray}\Lf=\Lf(G_n(p))
&=&\La(G_n(p))+p(\mathbf{J}_n-\mathbf{I}_n)\label{Equ-7}\\
&=&
\left(\mathbf{D}(G_n(p))-2\sum_{i=1}^n\sum_{j>i}a_{ij}/n\mathbf{I}_n\right)
-\big(\mathbf{A}(G_n(p))-p(\mathbf{J}_n-\mathbf{I}_n)\big).\nonumber
\end{eqnarray} The following result is concerned with the LSD of
$\Lf$.
\begin{thm}\label{LSD-LapEnergy}
Let $G_n(p)$ be a random graph of $\mathcal{G}_{n}(p)$. Then
$$\lim_{n\rightarrow\infty}\Phi_{(\sigma\sqrt{n})^{-1}\Lf}(x)
=\Psi(x)\mbox{ a.s.}$$
\end{thm}
To prove the theorem above, we introduce an auxiliary matrix as
follows:
\begin{eqnarray*}
\Ls=\Ls(G_n(p))
&=&\mathbf{L}(G_n(p))-(n-1)p\mathbf{I}_n+p(\mathbf{J}_n-\mathbf{I}_n)\\
&=& \big(\mathbf{D}(G_n(p))-(n-1)p\mathbf{I}_n\big)-
\big(\mathbf{A}(G_n(p))-p(\mathbf{J}_n-\mathbf{I}_n)\big).
\end{eqnarray*}

First of all, one can readily see that $\mathbf{L}_2$ is a Markov
matrix in which the Wigner matrix is
$-\mathbf{A}(G_n(p))+p(\mathbf{J}_n-\mathbf{I}_n)$ and the diagonal
matrix is $-\mathbf{D}(G_n(p))+(n-1)p\mathbf{I}_n$. Furthermore, the
off-diagonal entries of $\sigma^{-1}\mathbf{L}_2$ have mean 0 and
variance 1. Since the LSD $\Psi(x)$ does not depend on the random
variables $x_{ij}$, Theorem \ref{LSD-Markov} yields
$$\label{Equ-2}\lim_{n\rightarrow\infty}
\Phi_{(\sigma\sqrt{n})^{-1}\mathbf{L}_2}(x)=\Psi(x) \mbox{ a.s.}
$$

In what follows, we shall show that $(\sigma\sqrt{n})^{-1}\Lf$ and
$(\sigma\sqrt{n})^{-1}\Ls$ have the same LSD $\Psi(x)$, by which
Theorem \ref{LSD-LapEnergy} follows. To this end, we first estimate
the difference $(\sigma\sqrt{n})^{-1}(\mathbf{L}_1- \mathbf{L}_2) $
by Chernoff's inequality (see \cite{JLR}, pp. 26 for instance).

\begin{lem}[\bf Chernoff's Inequality]\label{Cher-Inequ}
Let $X$ be a random variable with binomial distribution $Bi(n,p)$.
Then, for any $\epsilon>0$,
$$
\Pro(|X-\E(X)|\ge
\epsilon)\le\exp\left\{-\frac{\epsilon^2}{2(np-\epsilon/3)}\right\}.$$
\end{lem}
Apparently,
$$
(\sigma\sqrt{n})^{-1}\mathbf{L}_2- (\sigma\sqrt{n})^{-1}\mathbf{L}_1
= (\sigma\sqrt{n})^{-1}
\left(2\sum_{i=1}^n\sum_{j>i}a_{ij}/n-(n-1)p\right)\mathbf{I}_n.
$$
Denote $(\sigma\sqrt{n})^{-1}
(2\sum_{i=1}^n\sum_{j>i}a_{ij}/n-(n-1)p)$ by $\Delta_n$ for
convenience. By means of Lemma \ref{Cher-Inequ}, for any given
$\epsilon>0$, we have
\begin{eqnarray*} & &\Pro\left((\sigma\sqrt{n})^{-1}
\left|2\sum_{i=1}^n\sum_{j>i}a_{ij}/n-(n-1)p\right|\ge\epsilon\right)\\
&=&\Pro\left(
\left|\sum_{i=1}^n\sum_{j>i}a_{ij}-\frac{n(n-1)p}{2}\right|\ge
\frac{\epsilon\cdot
\sigma n^{3/2}}{2}\right)\\
 &\le &
 \exp\left\{-\frac{2^{-2}(\epsilon\sigma)^2
 n^3}{2(n(n-1)p+\epsilon\sigma
n^{3/2}/6)}\right\}\\
 &<& \exp\left\{-\frac{(\epsilon\sigma)^2\cdot n^3}{8(p+\epsilon
\sigma/6)\cdot n^2}\right\}\\
&=& \exp\left\{-\frac{(\epsilon\sigma)^2}{8(p+\epsilon
\sigma/6)}\cdot n\right\}. \end{eqnarray*} Therefore, by the first
Borel-Cantelli lemma (see \cite{Bil}, pp. 59 for instance), we can
deduce
$$|\Delta_n|=(\sigma\sqrt{n})^{-1}
\left|2\sum_{i=1}^n\sum_{j>i}a_{ij}/n-(n-1)p\right|\rightarrow
0\mbox{ a.s. }(n\rightarrow\infty).
$$

Furthermore, it is easy to see that  $\lambda$ is an eigenvalue of
$(\sigma\sqrt{n})^{-1}\mathbf{L}_1$ if and only if
$\lambda+\Delta_n$ is an eigenvalue of
$(\sigma\sqrt{n})^{-1}\mathbf{L}_2$. By the definition of the ESD,
it follows that
\begin{equation}\label{Equ-1}
\Phi_{(\sigma\sqrt{n})^{-1}\mathbf{L}_1}(x)=
\Phi_{(\sigma\sqrt{n})^{-1}\mathbf{L}_2}(x+\Delta_n).
\end{equation}
Clearly, for any $\epsilon>0$, there exists $N$ such that
$|\Delta_n|<\epsilon \mbox{ a.s.}$ for all $n>N$. Noting that
$\Phi_{(\sigma\sqrt{n})^{-1}\mathbf{L}_2}(x)$ is an increasing
function, for all $n>N$, we have
$$\Phi_{(\sigma\sqrt{n})^{-1}\mathbf{L}_2}(x-\epsilon)\leq\Phi_{(\sigma\sqrt{n})^{-1}\mathbf{L}_2}(x+\Delta_n)
\leq\Phi_{(\sigma\sqrt{n})^{-1}\mathbf{L}_2}(x+\epsilon)\mbox{
a.s}.$$ Consequently,
\begin{eqnarray*}
\Psi(x-\epsilon)&=&\lim_{n\rightarrow\infty}\Phi_{(\sigma\sqrt{n})^{-1}\mathbf{L}_2}(x-\epsilon)\\
&\leq&\lim_{n\rightarrow\infty}\Phi_{(\sigma\sqrt{n})^{-1}\mathbf{L}_2}(x+\Delta_n)\\
&\leq&\lim_{n\rightarrow\infty}\Phi_{(\sigma\sqrt{n})^{-1}\mathbf{L}_2}(x+\epsilon)=\Psi(x+\epsilon)\mbox{
a.s}. \end{eqnarray*} Moreover, since the density of $\Psi(x)$ is
smooth bounded, $\Psi(x)$ is continuous. Together with the fact that
$\epsilon$ is arbitrary, we conclude
$$\lim_{n\rightarrow\infty}\Phi_{(\sigma\sqrt{n})^{-1}\mathbf{L}_1}(x)
=\lim_{n\rightarrow\infty}\Phi_{(\sigma\sqrt{n})^{-1}\mathbf{L}_2}(x+\Delta_n)=\Psi(x)
\mbox{ a.s}.,$$ which completes the proof of Theorem
\ref{LSD-LapEnergy}.

\subsection{The bounds of $\EnL(G_n(p))$}

In this part, we shall establish the lower and upper bounds of
$\EnL(G_n(p))$ by employing Theorem \ref{LSD-LapEnergy}, and then
show that Conjecture \ref{Gutman} is true for almost all graphs at
last.

Let $X$ be a random variable with the distribution $\Psi(x)$. We
start with an estimate of $\E|X|=\int |x|d\Psi(x)$. Since $\Psi(x)$
is the free convolution of the standard semi-circle distribution
$\Phi_{0,1}(x)$ and the standard normal measure, let us investigate
the free convolution in depth. Here, we follow the notation given by
Voiculescu \cite{V}. The Cauchy-Stieltjes transform of a probability
measure $\mu$ is
$$
G_{\mu}(z)=\int_{-\infty}^{\infty}\frac{\mu(dx)}{z-x}
$$
which is analytic on the complex upper half plane. For some
$\alpha,\beta>0$, there exists a domain
$D_{\alpha,\beta}=\{u+iv\mid|u|<\alpha v,v>\beta\}$ on which
$G_{\mu}$ is univalent. For the image $G_{\mu}(D_{\alpha,\beta})$,
we can define the inverse function $K_{\mu}$ of $G_{\mu}$ in the
area $\Gamma_{a,b}=\{u+iv\mid|u|<-av, -b<v<0\}$. And let
$R_{\mu}(z)=K_{\mu}(z)-1/z.$ Then for probability measures $\mu$ and
$\nu$, there exists a unique probability measure, denoted by
$\mu\boxplus\nu$, on $\Gamma_{a,b}$ such that
$$
R_{\mu\boxplus\nu}=R_{\mu}+R_{\nu}.
$$
The measure $\mu\boxplus\nu$ is said to be the free convolution of
$\mu$ and $\nu$.

In the above definition, the Cauchy-Stieltjes transform and inverse
function may be difficult to compute in practice. Consequently, we
do not compute $\E|X|$ directly. In what follows, we employ another
definition of free convolution via combinatorial way (see
\cite{wb,rs}) applicable only to probability measures with all
moments.

For  probability measure $\mu$, set $m_k=\int x^k\mu(dx)$ and
$$
M_\mu(z)=1+\sum_{k=1}^{\infty}m_kz^k.
$$
Define a formal power series
$$
T_\mu(z)=\sum_{k=1}^{\infty}c_kz^{k-1}
$$
such that
$$
M_\mu(z)=1+z M_\mu(z)T_\mu(zM_\mu(z)).
$$
Then, the free convolution of $\mu,\nu$ is the probability measure
$\mu\boxplus\nu$ satisfying
\begin{equation}\label{rt}
T_{\mu\boxplus\nu}(z)=T_{\mu}(z)+T_{\nu}(z).
\end{equation}
It is not difficult to see that this definition is coincident with
the analytical one (see \cite{rs}).

Next, we calculate $\E|X|$ by the following result due to Bryc
\cite{wb}. Let $M_{\mu,n}\equiv M_\mu(z)$ $\mbox{ mod }z^{n+1}$,
$T_{\mu,n}(z)\equiv T_\mu(z)\mbox{ mod }z^{n+1}$ be the $n$-th
truncations, {\it i.e.,} $M_{\mu,n}=1+\sum_{k=1}^{n}m_kz^k$ and
$T_{\mu,n}(z)=\sum_{k=1}^{n+1}c_kz^{k-1}$.
\begin{lem}[\bf Bryc \cite{wb}]
With $M_{\mu,0}(z)=1$ and $c_1=M'_{\mu,1}(0)$, we have
$$
M_{\mu,n}(z)\equiv1+zM_{\mu,n-1}(z)T_{\mu,n-1}(zM_{\mu,n-1}(z))\mod
z^{n+1},n\geq 1,
$$
and
$$
c_k=\left.-{\frac{1}{k-1}} {\frac{1}{k!}}{ \frac{d^k}{dz^k}}
\frac{1}{M_{\mu,n}^{k-1}(z)}\right|_{z=0}.
$$
\end{lem}

Therefore, combining with the formula (\ref{rt}), we can calculate
the moments of $\mu\boxplus\nu$ by the moments of $\mu$, $\nu$ in
recurrence. It is not difficult to verify that $\E X^2=2$ and $\E
X^4=9$ (see \cite{wb} for details). Employing Cauchy-Schwartz
inequality
$$
|\E(XY)|^2\leq \E X^2\cdot \E Y^2,
$$
we have $$\E |X|\leq \sqrt{\E X^2}$$ and
$$
(\E X^2)^2\leq \E |X|\cdot \E |X|^3\leq \E |X|\cdot \sqrt{\E
X^2\cdot \E X^4}.
$$
Therefore, $$\frac{2\sqrt{2}}{3}\le\E|X|\le\sqrt{2}.$$

In what follows, we shall establish the lower and upper bounds of
$\EnL(G_n(p))$ by employing an estimate of the energy $\En(\Lf)$. We
first investigate the convergence of $\int
|x|d\Phi_{(\sigma\sqrt{n})^{-1}\Lf}(x)$. Let $I$ be the interval
$[-1,1]$. By Theorem \ref{LSD-LapEnergy} and the bounded convergence
theorem, one can easily see that
\begin{equation}\label{Equ-8}
\lim_{n\rightarrow\infty}\int_I
|x|d\Phi_{(\sigma\sqrt{n})^{-1}\Lf}(x) =\int_I |x|d\Psi(x)\mbox{
a.s.}\end{equation} We proceed to prove that
$$\lim_{n\rightarrow\infty} \int_{I^c}
|x|d\Phi_{(\sigma\sqrt{n})^{-1}\Lf}(x)=\int_{I^c} |x|d\Psi(x)\mbox{
a.s.}$$where $I^c=\mathbb{R}\setminus I$. Since $\sigma^{-1}\Ls$ is
the Markov matrix such that the off-diagonal entries have mean 0 and
variance 1, we can deduce, by Equation (\ref{Moments}), that
\begin{equation}\label{Equ-20}
\lim_{n\rightarrow\infty}\int
x^2d\Phi_{(\sigma\sqrt{n})^{-1}\Ls}(x)=\int x^2d\Psi(x)\mbox{
a.s.}\end{equation} According to the relation (\ref{Equ-1}), we have
\begin{eqnarray*}
\int x^2d\Phi_{(\sigma\sqrt{n})^{-1}\Lf}(x)
&=&\int x^2d\Phi_{(\sigma\sqrt{n})^{-1}\Ls}(x+\Delta_n)\\
&=&\int
(x-\Delta_n)^2d\Phi_{(\sigma\sqrt{n})^{-1}\Ls}(x)\\
&=&\int x^2d\Phi_{(\sigma\sqrt{n})^{-1}\Ls}(x)-2\Delta_n\int
xd\Phi_{(\sigma\sqrt{n})^{-1}\Ls}(x)\\
&& +\Delta_n^2\int d\Phi_{(\sigma\sqrt{n})^{-1}\Ls}(x).
\end{eqnarray*} Since
$\lim_{n\rightarrow\infty}\Delta_n=0$ a.s., Equation (\ref{Equ-20})
implies that
\begin{equation}\label{Equ-3}\lim_{n\rightarrow\infty}\int
x^2d\Phi_{(\sigma\sqrt{n})^{-1}\Lf}(x)=\lim_{n\rightarrow\infty}\int
x^2d\Phi_{(\sigma\sqrt{n})^{-1}\Ls}(x)=\int x^2d\Psi(x)\mbox{
a.s.}\end{equation} Consequently,
\begin{eqnarray*}\lim_{n\rightarrow\infty}\int_{I^c}
x^2d\Phi_{(\sigma\sqrt{n})^{-1}\Lf}(x)&=&\lim_{n\rightarrow\infty}\left(\int
x^2d\Phi_{(\sigma\sqrt{n})^{-1}\Lf}(x)-\int_{I}
x^2d\Phi_{(\sigma\sqrt{n})^{-1}\Lf}(x)\right)\nonumber\\
&=&\int_{I^c} x^2d\Psi(x)\mbox{ a.s.}\label{Equ-9}\end{eqnarray*}
\begin{lem}[\bf Billingsley \cite{Bil} pp. 219]\label{Billinsley}
Let $\mu$ be a measure. Suppose that functions $a_n, b_n, f_n$
converges almost everywhere to functions $a, b, f$, respectively,
and that $a_n\le f_n\le b_n$ almost everywhere. If $\int
a_nd\mu\rightarrow\int a\hspace{2pt}d\mu$ and $\int
b_nd\mu\rightarrow\int b\hspace{2pt}d\mu$, then $\int
f_nd\mu\rightarrow\int fd\mu$.\end{lem} Suppose that
$\phi_{(\sigma\sqrt{n})^{-1}\Lf}(x)$ is the density of
$\Phi_{(\sigma\sqrt{n})^{-1}\Lf}(x)$. By virtue of Theorem
\ref{LSD-LapEnergy} and Lemma \ref{Billinsley}, we can deduce by
setting $a_n(x)=0$, $b_n(x)=x^2\phi_{(\sigma\sqrt{n})^{-1}\Lf}(x)$
and $f_n(x)=|x|\phi_{(\sigma\sqrt{n})^{-1}\Lf}(x)$ that
$$\lim_{n\rightarrow\infty} \int_{I^c}
|x|d\Phi_{(\sigma\sqrt{n})^{-1}\Lf}(x)=\int_{I^c} |x|d\Psi(x)\mbox{
a.s.}$$ Combining the above equation with Equation (\ref{Equ-8}), we
have$$\lim_{n\rightarrow\infty} \int
|x|d\Phi_{(\sigma\sqrt{n})^{-1}\Lf}(x)=\int |x|d\Psi(x)\mbox{
a.s.}$$

We are now ready to present an estimate of the energy $\En(\Lf)$. By
an argument similar to evaluate the energy
$\En(\overline{\mathbf{A}})$, we have
\begin{eqnarray*}\En(\Lf)/\sigma n^{3/2} &=&
\int |x|d\Phi_{(\sigma\sqrt{n})^{-1}\Lf}(x)\\
&\rightarrow& \int |x|d\Psi(x)\mbox{ a.s. }(n\rightarrow\infty).\\
\end{eqnarray*} Since $2\sqrt{2}/3\le \E|X|\le \sqrt{2}$,
$$
\frac{2\sqrt{2}}{3}\le\frac{\En(\Lf)}{\sigma n^{3/2}}\le\sqrt{2}
\mbox{ a.s. }(n\rightarrow\infty).$$ Consequently,
\begin{equation}\label{Equ-N}
\left(\frac{2\sqrt{2}}{3}\sigma+o(1)\right)n^{3/2}\le\En(\Lf)
\le\left(\sqrt{2}\sigma+o(1)\right)n^{3/2} \mbox{ a.s.
}\end{equation}

Employing the equation above, we can establish the lower and upper
bounds of $\EnL(G_n(p))$. Note that $\EnL(G_n(p))=\En(\La)$
according to the definition of the energy of a matrix. So we turn
our attention to the bounds of $\En(\La)$. By means of Equation
(\ref{Equ-7}), we have
$$\Lf=\La+p(\mathbf{J}_n-\mathbf{I}_n)\mbox{ and }
\La=\Lf+p(\mathbf{I}_n-\mathbf{J}_n).$$ Thus, Lemma \ref{ky} yields
that
$$\En(\Lf)-\En(p(\mathbf{J}_n-\mathbf{I}_n))\le
\En(\La)\le\En(\Lf)+\En(p(\mathbf{I}_n-\mathbf{J}_n)).$$ Recalling
the fact that $\En(p(\mathbf{J}_n-\mathbf{I}_n))=
\En(p(\mathbf{I}_n-\mathbf{J}_n))=2p(n-1)$, Equation (\ref{Equ-N})
implies that
$$\left(\frac{2\sqrt{2}}{3}\sigma +o(1)\right)n^{3/2}-2p(n-1)\le\En(\La)
\le\left(\sqrt{2}\sigma+o(1)\right)n^{3/2}+2p(n-1)\mbox{ a.s. }$$
Therefore, we obtain the the lower and upper bounds of the Laplacian
energy for almost all graphs.
\begin{thm}
Almost every random graph $G_n(p)$ satisfies
$$
\left(\frac{2\sqrt 2}{3}\sigma+o(1)\right)\cdot n^{3/2}\leq
\EnL(G_n(p))\leq \left(\sqrt{2}\sigma+o(1)\right)\cdot n^{3/2}.
$$
\end{thm}

Since a.e. random graph $G_n(p)$ satisfies
$$
\lim_{n\rightarrow\infty}\frac{\En(G_n(p))}{n^{3/2}}
=\frac{8}{3\pi}\sigma <\frac{2\sqrt 2}{3}\sigma
\leq\lim_{n\rightarrow\infty}\frac{\EnL(G_n(p))}{n^{3/2}},
$$
we thus establish the result below.
\begin{thm}
For almost every random graph $G_n(p)$, $\En(G_n(p))<\EnL(G_n(p))$.
\end{thm}

By virtue of the theorem above, Conjecture \ref{Gutman} is true
for almost all graphs.\\

\noindent {\bf Acknowledgement:} The authors are very grateful to
the referees for detailed suggestions and comments, which helped to
improve the presentation of the manuscript significantly.

\end{document}